\documentclass[12pt]{article}
\usepackage{enumerate}
\usepackage{amsfonts}
\usepackage{amsmath}
\usepackage{amssymb}
\usepackage{amsthm}

\def\D{\mathcal{D}}
\def\V{\mathcal{V}}

\def\H{\mathcal{H}}
\def\K{\mathcal{K}}

\def\T{\mathfrak{T}}
\def\B{\mathfrak{B}}

\newcommand{\rank}{\mathrm{rank}}

\newcommand{\Tr}{\mathrm{Tr}}

\newcommand{\shs}{\hspace{1pt}}

\newcounter{defin}  \newcounter{lemma}  \newcounter{theorem}
\newcounter{property} \newcounter{corol}  \newcounter{remark} \newcounter{example}

\newenvironment{lemma}{\par\refstepcounter{lemma}
     \textbf{Lemma \thelemma.} }{\rm\par}

\newenvironment{property}{\par\refstepcounter{property}
     \textbf{Proposition \theproperty.}\ }{\rm\par}
\newenvironment{corollary}{\par\refstepcounter{corol}
     \textbf{Corollary \thecorol.} }{\rm\par}

\begin{document}

\title{Relatively bounded operators and the operator E-norms (addition to arXiv:1806.05668)}

\author{M.E. Shirokov\footnote{Steklov Mathematical Institute, RAS, Moscow, email:msh@mi.ras.ru}}
\date{}
\maketitle

\vspace{-10pt}

\begin{abstract}
In this brief note we describe relations between the well known notion of a relatively bounded operator and the operator
E-norms considered in [arXiv:1806.05668].

We show that the set of all $\sqrt{G}$-bounded operators equipped with the E-norm induced by a positive operator $G$
is the Banach space of all operators with finite E-norm and that the $\sqrt{G}$-bound is a continuous seminorm on this space.

We also show that the set of all $\sqrt{G}$-infinitesimal operators (operators with zero $\sqrt{G}$-bound) equipped with the E-norm induced by a positive operator $G$ is the completion of the algebra $\B(\H)$ of bounded operators w.r.t. this norm.

Some properties of $\sqrt{G}$-infinitesimal operators are considered.
\end{abstract}

\section{Introduction}  A linear operator $A$ on a Hilbert space $\H$ is called relatively bounded w.r.t. a linear operator $B$ (briefly, $B$-bounded) if
$\mathcal{D}(B)\subseteq \D(A)$ and
\begin{equation}\label{rb-rel-}
\|A\varphi\|^2\leq a^2\|\varphi\|^2+b^2\|B\varphi\|^2\quad \forall \varphi\in\mathcal{D}(B)
\end{equation}
for some  nonnegative numbers $a$ and $b$. The infimum over $b$ for which (\ref{rb-rel-}) holds (with some $a$) is called the $B$-bound for $A$.
If this $B$-bound is equal to zero then $A$ is called $B$-infinitesimal operator (infinitesimally bounded w.r.t. $B$). These notions are widely used in
the modern operator theory, in particular, in analysis of perturbations of unbounded operators in a Hilbert space \cite{Kato,BS}.

The operator E-norm $\|A\|_E^G$ of a bounded operator $A$ on a Hilbert space $\H$ induced by a positive unbounded operator $G$ is introduced\footnote{As far as I know. I would be grateful for any references.} in \cite{CSR} as the maximum of $\sqrt{\Tr A \rho A^*}$ over all states (positive operators with unit trace) $\rho$ such that $\Tr G\rho\leq E$.\footnote{The value of $\Tr {G}\rho$ (finite or infinite) is defined as $\sup_n\Tr P_nG\rho$, where $P_n$ is the spectral projector of $G$ corresponding to the interval $[0,n]$.} This norm was used in \cite{CSR} to obtain the modification of the Kretschmann-Schlingemann-Werner theorem\footnote{The original Kretschmann-Schlingemann-Werner theorem obtained in \cite{Kr&W} quantifies  continuity of the Stinespring dilation of CP linear maps w.r.t. the diamond norm ($cb$-norm) topology on the set of CP linear maps and the operator norm topology on the set of Stinespring operators.} which quantifies continuity of the Stinespring dilation of a quantum channel w.r.t. the strong convergence topology on the set of  channels and the strong operator topology on the set of Stinepring isometries.\footnote{If $G$ is a unbounded operator with discrete spectrum of finite multiplicity then any of the E-norms $\|\cdot\|_E^G$, $E>0$, generates the strong operator topology on the unit ball of $\B(\H)$ \cite[Proposition 2]{CSR}.} These norms
are studied in detail in \cite{ECN}, where they are extended to unbounded operators. The extended operator E-norms
and the corresponding Banach spaces of unbounded operators have different applications described in \cite[Section 5]{ECN},\cite{SPM}.

In this note we describe relations between the notion of a relatively bounded operator and the operator
E-norms extended to unbounded operators.

We show that the set of all $\sqrt{G}$-bounded operators equipped with the E-norm induced by a positive operator $G$
coincides with the Banach space of all operators with finite E-norm denoted by $\B_G(\H)$ in \cite{ECN} and that the $\sqrt{G}$-bound of an operator is a continuous seminorm on $\B_G(\H)$. We obtain an explicit formula for the E-norm $\|A\|_E^G$ in terms of the set of coefficients $(a,b)$ for which (\ref{rb-rel-}) holds with $B=\sqrt{G}$ and the expression for the $\sqrt{G}$-bound of an operator $A$ via  $\|A\|_E^G$.

We also show that the set of all $\sqrt{G}$-infinitesimal operators (operators with zero $\sqrt{G}$-bound) equipped with the E-norm $\|\cdot\|_E^G$ coincides with the completion of the algebra $\B(\H)$ of bounded operators w.r.t. this norm denoted by $\B^0_G(\H)$ in \cite{ECN}.

\section{Definitions and the main result}

Let $\mathcal{H}$ be a separable infinite-dimensional Hilbert space, $\mathfrak{B}(\mathcal{H})$
-- the algebra of all bounded operators on $\mathcal{H}$ with the operator norm $\|\!\cdot\!\|$ and $\mathfrak{T}(\mathcal{H})$ --
the Banach space of all trace-class operators on $\mathcal{H}$ with the trace norm $\|\!\cdot\!\|_1$ (the Schatten class of order 1) \cite{Kato,BS}. Let
$\mathfrak{S}(\mathcal{H})$ be the set of quantum states -- positive operators in
$\mathfrak{T}(\mathcal{H})$ with unit trace \cite{H-SCI}.

Let ${G}$ be a positive (semidefinite) operator on $\H$ with a dense domain $\mathcal{D}({G})$ such that
\begin{equation}\label{G-cond}
\inf\left\{\shs\|G\varphi\|\,|\,\varphi\in\mathcal{D}({G}),\|\varphi\|=1\shs\right\}=0.
\end{equation}
For a given linear (bounded or unbounded) operator $A$ such that $\mathcal{D}(\sqrt{G})\subseteq \D(A)$ the operator E-norm induced by $G$
is defined in \cite{ECN} as
\begin{equation}\label{ec-on}
 \|A\|^{G}_E\doteq \sup\left.\left\{\sqrt{\Tr A\rho A^*} \,\right|\, \rho\in\mathfrak{S}(\mathcal{H}),\Tr {G}\rho\leq E, \rank \rho<+\infty\right\}
\end{equation}
where we assume that\footnote{This assumption is made to avoid the notion of adjoint operator.}
\begin{equation}\label{ab-d}
  A\rho A^*\doteq\sum_i|\alpha_i\rangle\langle\alpha_i|,\quad\quad |\alpha_i\rangle=A|\varphi_i\rangle,
\end{equation}
provided that $\rho=\sum_i |\varphi_i\rangle\langle\varphi_i|$ (by using Schrodinger's mixture theorem (see \cite[Ch.8]{B&Z}) it is easy to show that the r.h.s. of (\ref{ab-d}) does not depend on this decomposition of $\rho$).\footnote{For any vector $|\alpha\rangle$ the symbol $|\alpha\rangle\langle\alpha|$ denotes the 1-rank operator mapping a vector $|\beta\rangle$ to $\langle\alpha|\beta\rangle|\alpha\rangle$.} By using purification of a state it is easy to see that
\begin{equation}\label{ec-on+}
 \|A\|^{G}_E=\sup_n\sup\left\{\|A\otimes I_{\H_n}\varphi\|\,\left|\, \varphi\in\H\otimes\H_n, \|\varphi\|=1, \|\sqrt{G}\otimes I_{\H_n}\varphi\|\leq \sqrt{E}\,\right\}\right.,
\end{equation}
where $\H_n$ and $I_{\H_n}$ denote, respectively, a $n$-dimensional
Hilbert space and the unit operator in this space.

For any given operator $A$ the nonnegative nondecreasing function $E\mapsto\left[\|A\|^{G}_E\right]^2$ is concave on $\mathbb{R}_+$ and tends to $\|A\|\leq+\infty$ as $E\rightarrow+\infty$ \cite{ECN}. This implies  that
\begin{equation}\label{E-n-eq}
\|A\|^{G}_{E_1}\leq \|A\|^{G}_{E_2}\leq \sqrt{E_2/E_1}\|A\|^{G}_{E_1}\quad\textrm{ for any } E_2>E_1>0.
\end{equation}
So, for given $G$ all the norms $\|A\|^{G}_{E}$ are equivalent. In particular, if $\|A\|^{G}_{E}$ is finite
for some $E>0$ then $\|A\|^{G}_{E}$ is finite for all $E>0$.

It is shown in \cite{ECN} that the set of all operators $A$ with finite $\|A\|^{G}_{E}$ equipped with the norm $\|\cdot\|^{G}_{E}$ and naturally defined linear operations is a nonseparable Banach space denoted therein by
$\B_G(\H)$.\footnote{We identify operators coinciding on the set $\D(\sqrt{G})$.} The completion $\B^0_G(\H)$ of the algebra $\B(\H)$ w.t.r. to the norm $\|\cdot\|^{G}_{E}$ is a proper subspace of $\B_G(\H)$ determined by
the condition
\begin{equation}\label{s-cond}
 \|A\|^{G}_E=o\shs(\sqrt{E})\quad\textup{ as }\quad E\rightarrow+\infty.
\end{equation}
If $G$ is a unbounded operator with discrete spectrum of finite multiplicity then the Banach space $\B^0_G(\H)$ is separable and for any $A\in\B^0_G(\H)$ its  E-norm can be defined by the simple formula
\begin{equation}\label{ec-on++}
 \|A\|^{G}_E=\sup\left\{\|A\varphi\|\,\left|\, \varphi\in\H, \|\varphi\|=1, \|\sqrt{G}\varphi\|\leq \sqrt{E}\,\right\}\right.,
\end{equation}
which means that the first supremum in (\ref{ec-on+}) is achieved at $n=1$ and the supremum in (\ref{ec-on}) can be taken over pure states \cite[Theorem 3F]{ECN}. Due to the assumption (\ref{G-cond}) the condition
$\|\varphi\|=1$ in (\ref{ec-on+}) and (\ref{ec-on++}) can be replaced by $\|\varphi\|\leq1$ \cite[Proposition 3A]{ECN}.\footnote{The question about  coincidence of (\ref{ec-on+}) and (\ref{ec-on++}) for any positive operator $G$ is open. It is easy to show that this coincidence is equivalent to concavity of the r.h.s. of (\ref{ec-on++}) as a function of $E$.}

According to the general definition mentioned in the Introduction an operator $A$ is called relatively bounded w.r.t. the operator $\sqrt{G}$ (briefly, $\sqrt{G}$-bounded) if
$\mathcal{D}(\sqrt{G})\subseteq \D(A)$ and
\begin{equation}\label{rb-rel}
\|A\varphi\|^2\leq a^2\|\varphi\|^2+b^2\|\sqrt{G}\varphi\|^2,\quad \forall \varphi\in\mathcal{D}(\sqrt{G})
\end{equation}
for some  nonnegative numbers $a$ and $b$. Denote by $\Gamma_{\!\sqrt{G}}(A)$ the set of all pairs $(a,b)$ for which
(\ref{rb-rel}) holds. It is easy to see that $\Gamma_{\!\sqrt{G}}(A)$ is a closed subset of $\mathbb{\mathbb{R}}^2_+$. The  $\sqrt{G}$-bound of $A$ (denoted by $b_{\sqrt{G}}(A)$ in what follows) is defined as
$$
b_{\sqrt{G}}(A)=\inf\left\{b\,|\,(a,b)\in\Gamma_{\!\sqrt{G}}(A)\right\}.
$$

\begin{lemma}\label{bl}
\emph{A pair $(a,b)$ belongs to the set $\,\Gamma_{\!\sqrt{G}}(A)$ if and only if $\|A\|^{G}_{E}\leq \sqrt{a^2+b^2 E}$ for all $\,E>0$.}
\end{lemma}\smallskip

\emph{Proof.} If $\|A\|^{G}_{E}\leq \sqrt{a^2+b^2 E}$ then definition (\ref{ec-on+}) implies that
$$
\|A\varphi\|\leq \|A\|^{G}_{\|\sqrt{G}\varphi\|^2}\leq \sqrt{a^2+b^2 \|\sqrt{G}\varphi\|^2}
$$
for any unit vector $\varphi$ in $\D(\sqrt{G})$. Hence $(a,b)\in\Gamma_{\!\sqrt{G}}(A)$.

If $(a,b)\in\Gamma_{\!\sqrt{G}}(A)$ then it is easy to show that $(a,b)\in\Gamma_{\!\sqrt{G}\otimes I_{\H_n}}(A\otimes I_{\H_n})$, where $\H_n$ is a $n$-dimensional
Hilbert space, for any $n$ \cite[Theorem 7.1.20]{BS}. Hence
$$
\sup\left\{\|A\otimes I_{\H_n}\varphi\|\,\left|\,  \varphi\in\H\otimes\H_n, \|\varphi\|=1, \|\sqrt{G}\otimes I_{\H_n}\varphi\|\leq \sqrt{E}\right\}\right.\leq\sqrt{a^2+b^2 E}
$$
for any $n$ and $E>0$. So, definition (\ref{ec-on+}) implies that $\|A\|^{G}_{E}\leq \sqrt{a^2+b^2 E}$. $\square$\smallskip

\begin{property}\label{main-p} A) \emph{The Banach space $\B_G(\H)$ coincides (as a set) with the set of all $\sqrt{G}$-bounded operators. If $A\in\B_G(\H)$ and $E>0$ then}
$$
\|A\|^{G}_{E}=\inf\left.\left\{\sqrt{a^2+b^2 E}\;\right|(a,b)\in\Gamma_{\!\sqrt{G}}(A)\right\}\quad and \quad b_{\sqrt{G}}(A)=\lim_{E\rightarrow+\infty}\|A\|^{G}_{E}/\sqrt{E}.
$$
\emph{The limit in the last formula can be replaced by the infimum over all $\,E>0$.} \smallskip

B) \emph{The completion $\B^0_G(\H)$ of $\B(\H)$ w.r.t. the norm $\|\cdot\|_E^G$ coincides (as a set) with the set of all $\sqrt{G}$-infinitesimal operators, i.e. operators with
the $\sqrt{G}$-bound equal to $0$.}\smallskip

C) \emph{The function $b_{\sqrt{G}}(\cdot)$ is a continuous seminorm on $\B_G(\H)$ s.t. $b_{\sqrt{G}}^{-1}(0)=\B^0_G(\H)$. Quantitatively,
\begin{equation}\label{b-cb}
\left|\shs b_{\sqrt{G}}(A)-b_{\sqrt{G}}(B)\right|\leq b_{\sqrt{G}}(A-B)\leq \|A-B\|^{G}_{E}/\sqrt{E}
\end{equation}
for arbitrary $A,B$ in $\B_G(\H)$ and any $E>0$.}
\end{property}\medskip

\emph{Proof.} Since $E\mapsto [\|A\|^{G}_{E}]^2$ is a concave nonnegative function on $\mathbb{R}_+$, it coincides with the infimum of all linear functions $E\mapsto a^2+b^2 E$
such that $[\|A\|^{G}_{E}]^2\leq a^2+b^2 E$ for all $E>0$ and the function $E\mapsto [\|A\|^{G}_{E}]^2/E$ is non-increasing. So, the assertions A and B can be easily derived from Lemma \ref{bl}.

To prove C note first that the seminorm properites of $b_{\sqrt{G}}(\cdot)$ follow from the second formula in A, while B implies $\,b_{\sqrt{G}}^{-1}(0)=\B^0_G(\H)$. So, since the function $E\mapsto [\|A\|^{G}_{E}]^2/E$ is non-increasing for any given $A$, the inequality (\ref{b-cb}) follows from the triangle inequality for $b_{\sqrt{G}}(\cdot)$. $\square$\smallskip

Due to Proposition \ref{main-p} one can reformulate the results in Section 4 in \cite{ECN} using the notions  of      
$\sqrt{G}$-bounded and $\sqrt{G}$-infinitesimal operators. In particular, Theorem 3 in \cite{ECN} implies the following characterization of $\sqrt{G}$-infinitesimal operators.\smallskip
\begin{corollary}\label{main-c}
\emph{An operator $A$ defined on $\D(\sqrt{G})$ is $\sqrt{G}$-infinitesimal if and only if for
any separable Hilbert space $\K$  the operator $A\otimes I_{\K}$ (naturally defined on the set $\,\D(\sqrt{G})\otimes\K$)  has a continuous linear extension to the set 
$$
\V_{\sqrt{G}\otimes I_{\K},E}=\left\{\shs\varphi\in\H\otimes\K\,\left|\,\|\sqrt{G}\otimes I_{\K} \varphi\|^2\leq E\shs\right.\right\}
$$ for any $E>0$. If $A$ is a $\sqrt{G}$-infinitesimal operator then
\begin{equation}\label{tin}
\|A\otimes I_{\K}(\varphi-\psi)\|\leq \varepsilon\|A\|^{G}_{4E/\varepsilon^2}
\end{equation}
for any $\,\varphi$ and $\,\psi$ in $\V_{{G}\otimes I_{\K},E}$ such that $\|\varphi-\psi\|\leq\varepsilon$. The r.h.s. of (\ref{tin}) tends to zero as $\,\varepsilon\rightarrow 0^+$ by condition (\ref{s-cond}).}
\end{corollary}\smallskip

Proposition 6 in \cite{ECN} and Proposition \ref{main-p}B imply that any 2-positive linear map $\Phi:\B(\H)\rightarrow\B(\H)$ such that $\,\Phi(I_{\H})\leq I_{\H}\,$ having the predual map $\Phi_*:\T(\H)\rightarrow\T(\H)$
with finite\footnote{It is easy to see that $E\mapsto Y_{\Phi}(E)$ is a concave function. So, finiteness of $Y_{\Phi}(E)$ for some $E>0$ implies finiteness of $Y_{\Phi}(E)$ for all $E>0$ and boundedness of the function $E\mapsto Y_{\Phi}(E)/E$.}
\begin{equation*}
Y_{\Phi}(E)\doteq \sup\left\{\shs\Tr {G}\Phi_*(\rho)\,|\,\rho\in\mathfrak{S}(\mathcal{H}),\Tr {G}\rho\leq E\,\right\}
\end{equation*}
is uniquely extended to a linear transformation of the set of all $\sqrt{G}$-infinitesimal operators bounded w.r.t.
the norm $\|\cdot\|^{G}_E$.\footnote{If $G$ is a Hamiltonian of a quantum system then the quantity $Y_{\Phi}(E)/E$ can be treated as an energy amplification factor of $\Phi$. Quantum channels $\Phi$ with finite $Y_{\Phi}(E)$ called energy-limited in \cite{W-EBN} naturally appear as realistic quantum dynamical maps.}\smallskip

Finally, consider application of the formula for the $\sqrt{G}$-bound  in Proposition \ref{main-p}A.  

\textbf{Example.} Let $\H = L_2(\mathbb{R})$  and $S(\mathbb{R})$ be the set of infinitely differentiable rapidly decreasing functions
with all the derivatives tending to zero quicker than any degree of $|x|$ when $|x|\rightarrow+\infty$.
Consider the operators $q$ and $p$ defined on the set $S(\mathbb{R})$ by setting
$$
(q\varphi)(x) = x\varphi(x)\quad\textrm{and}\quad  (p\shs\varphi)(x) = \frac{1}{i}\frac{d}{dx}\varphi(x).
$$
These operators are essentially self-adjoint. They represent (sharp) real
observables of position and momentum of a quantum particle in the system of units where Planck's
constant $\hbar$ is equal to $1$ \cite[Ch.12]{H-SCI}. For given $\omega>0$ consider the operators
\begin{equation}\label{a-oper-d}
a=(\omega q+ip)/\sqrt{2\omega}\quad \textrm{and} \quad a^{\dagger}=(\omega q-ip)/\sqrt{2\omega}
\end{equation}
defined on $S(\mathbb{R})$. The operator $N=a^{\dagger}a=aa^{\dagger}-I_{\H}$ 
is positive and essentially self-adjoint. It represents (sharp) real
observable of the number of quanta of the harmonic oscillator with frequency $\omega$.
In \cite[Section 5]{ECN} the following estimates are obtained
\begin{equation}\label{QP-en}
\sqrt{\frac{2E+1/2}{\omega}}< \|q\|^{N}_E\leq\sqrt{\frac{2E+1}{\omega}},\quad \sqrt{(2E+1/2)\omega}< \|p\|^N_E\leq\sqrt{(2E+1)\omega}
\end{equation}
(the \emph{E}-norms of $q$ and $p$ depend on $\omega$, since the operator $N$ depends on $\omega$). Thus, 
the second formula in Proposition \ref{main-p}A implies that
$b_{\sqrt{N}}(p)=\sqrt{2/\omega}$ and $b_{\sqrt{N}}(q)=\sqrt{2\omega}$. $\square$

\bigskip

I am grateful to T.V.Shulman for the help and useful discussion.

\end{document}